\input amstex 
\documentstyle{amsppt}
\magnification=\magstep1
\pageheight{23truecm}\pagewidth{16truecm}
\parindent=10pt
 2 
 1 
 1 
 
\define\oldqed{\hbox{\vrule height7pt width3pt depth1pt}}
\def\fra{\frak a}
\def\frm{\frak m}
\def\frq{\frak q}
\def\frn{\frak n}

\def\FA{\bold F_{\!A}}
\def\FAe{\bold F_{\!A}^e}
\def\FR{\bold F_{\!R}}
\def\FS{\bold F_{\!S}}
\def\fraM{\frak M}

\def\bbP{\Bbb P}
\def\bbQ{\Bbb Q}
\def\bbZ{\Bbb Z}
\def\1A{\!\,^1\!A}
\def\eA{\!\,^e\!A}

\def\ehk{ e_{\operatorname{HK}}}
\def\mhk{\operatorname{m_{\operatorname{HK}}}}
\def\relmhk{\operatorname{rel.m_{\operatorname{HK}}}}
\def\ann{\operatorname{ann}}
\def\Ass{\operatorname{Ass}}

\def\chara{\operatorname{char}}

\def\Cl{\operatorname{Cl}}
\def\cl{\operatorname{cl}}

\def\divi{\operatorname{div}}

\def\Hom{\operatorname{Hom}}

\def\Ker{\operatorname{Ker}}
\def\Min{\operatorname{Min}}
\def\ord{\operatorname{ord}}
\def\pd{\operatorname{pd}}

\def\rank{\operatorname{rank}}
\def\Soc{\operatorname{Soc}}

\def\Supp{\operatorname{Supp}}
\def\Sym{\operatorname{Sym}}
\topmatter
\title 
Minimal Hilbert-Kunz multiplicity 
\endtitle
\date
Mar. 7, 2003
\enddate 
\rightheadtext{Minimal Hilbert-Kunz multiplicity }
\affil
Kei-ichi Watanabe \\
Department of Mathematics \\
College of Humanities and Sciences, 
Nihon University \\
Setagaya-ku, Tokyo 156--0045, Japan \\
e-mail: {\tt watanabe\@math.chs.nihon-u.ac.jp}
\\ \\ and  \\ \\
Ken-ichi Yoshida \\
Graduate School of Mathematics, Nagoya University \\
Chikusa-ku, Nagoya  464--8602, Japan \\
e-mail: {\tt yoshida\@math.nagoya-u.ac.jp}
\endaffil
\abstract
In this paper, we ask the following question: 
what is the minimal value of the difference 
$\ehk(I) - \ehk(I')$ for ideals $I' \supseteq I$ with $l_A(I'/I) =1$? 
\par 
In order to answer to this question, 
we define the notion of {\it minimal Hilbert-Kunz multiplicity}  
for strongly F-regular rings.  
Moreover, we calculate this invariant  
for quotient singularities and for the coordinate ring 
of the Segre embedding: $\bbP^{r-1} \times \bbP^{s-1} \hookrightarrow
 \bbP^{rs-1}$, respectively.  
\endabstract
\endtopmatter
\document
\head 
{\bf Introduction}
\endhead 
Throughout this paper, let $A$ be a Noetherian ring of positive 
prime characteristic $p$. 
The purpose of this paper is to introduce the notion of minimal 
Hilbert-Kunz multiplicity which is a new invariant of local rings in positive 
characteristic. 
\par
The notion of Hilbert-Kunz multiplicity has been introduced by 
Kunz \cite{Ku1} in 1969, and has been studied in detail 
by Monsky \cite{Mo}; 
see also e.g. \cite{BC}, \cite{BCP}, \cite{Co}, \cite{HM}, \cite{Se1}, 
or \cite{WY1,WY2,WY3}. 
\par
Further, Hochster and Huneke \cite{Hu2} have pointed out that the tight closure 
$I^{*}$ of $I$ is the largest ideal containing $I$ having the same 
Hilbert-Kunz multiplicity as $I$; see Lemma 1.3. 
Thus it seems to be important to understand the 
Hilbert-Kunz multiplicity well. 
For example, the authors \cite{WY1} have proved that an unmixed local ring whose 
Hilbert-Kunz multiplicity one is regular. 
Also, they \cite{WY3} have given a formula of $\ehk(I)$ 
for any integrally closed ideal $I$ in a two-dimensional F-rational 
double point using McKay correspondence, Riemann--Roch formula. 
\par
One of the most important conjectures about Hilbert-Kunz multiplicities 
is that it is always a rational number. 
Let $A$ be a local ring and $I,\,J$ be $\frm$-primary ideals in $A$. 
Also, suppose that $J$ is a parameter ideal. 
Then it is known that $\ehk(J) = e(J)$, the usual multiplicity (and hence 
$\ehk(J)$ is an integer). 
In order to investigate the value of $\ehk(I)$, we study 
the difference \lq\lq$\ehk(J) - \ehk(I)$''. 
So it is natural to ask the following question. 

\proclaim {Question} 
What is the minimal value of the difference $\ehk(I)-\ehk(I')$ 
for ideals $I'\supseteq I$ with $l_A(I'/I) =1$.  
\endproclaim 
\par
To give an answer to this question, we introduce the notion of 
{\it relative minimal Hilbert-Kunz multiplicity} $\relmhk(A)$ and 
that of {\it minimal Hilbert-Kunz multiplcity} $\mhk(A)$ 
as follows:
$$
\aligned 
\relmhk(A) &:= \inf\{\ehk(I) -\ehk(I') \,\:\, I \subseteq I' \;\text{with}\; 
l_A(I'/I) =1\} \cr 
\mhk(A) &:= \liminf_{e \to \infty} \dfrac{l_A(A/\ann_A F_A^e(z))}{p^{ed}}, 
\endaligned 
$$ 
where $z$ is a generator of the socle of the 
injective hull $E_A(A/\frm)$. 
\newline
Note that $\relmhk(A) \ge \mhk(A)$; see Proposition 1.10.
\par
In Section 2, we prove that if $A$ is a Gorenstein local ring then 
$$
\ehk(J) -\ehk(J:\frm) = \relmhk(A) = \mhk(A)
$$ 
for any parameter ideal $J$ of $A$; see Theorem 2.1 for details. 
Note that a similar result is independently proved by Huneke and 
Leuschke \cite{HuL}.  
\par 
In general, if $A$ is not weakly F-regular (see Definition 1.1), then 
$\relmhk(A)=0$ (and hence $\mhk(A) =0$). 
Thus it suffices to consider weakly F-regular local rings in our context. 
\par
In Section 3, we will give a formula for minimal Hilbert-Kunz 
multiplicities of the canonical cover of $\bbQ$-Gorenstein F-regular local 
rings as follows$:$ 

\proclaim {Theorem 3.1} 
Let $A$ be a $\bbQ$-Gorenstein strongly F-regular local ring. 
Also, let $B = A \oplus K_At \oplus K_A^{(2)}t^2 \oplus 
\cdots \oplus K_A^{(r-1)}t^{r-1}$, the canonical cover of $A$, 
where $r = \ord(cl(K_A))$, $K_A^{(r)} = fA$ and $ft^r =1$. 
Also, suppose that $(r,p)= 1$. 
Then we have 
$$
\mhk(B) = r \cdot \mhk(A).
$$ 
\endproclaim 
\par
In Section 4, as an application of Theorem 3.1, we will give a formula 
for the minimal Hilbert-Kunz multiplicities of quotient singularities. 

\proclaim {Theorem 4.2} 
Let $k$ be a field of characteristic $p >0$, 
and let $A = k[x_1,\,\ldots,x_d]^G$ be the invariant subring by a 
finite subgroup $G$ of $GL(d,k)$ with $(p,\,|G|)=1$. 
Also, assume  that $G$ contains no pseudo-reflections.  
Then $\mhk(A) = 1/|G|$.  
\endproclaim 
\par
In Section 5, we will give a formula for the minimal Hilbert-Kunz multiplicities
of Segre products. 

\proclaim {Theorem 5.6}
Let $A=k[x_1,\,\ldots,x_r]\#k[y_1,\,\ldots,y_s]$, where 
$2 \le r \le s$, and put $d = r+s-1$. 
Then 
$$ 
\mhk(A)  =\frac{r!}{d!}\,S(d,r) 
+\frac{1}{d!}\sum_{k=1}^{r-1} \sum_{j=1}^{r-k}
\binom{r}{k+j} \binom{s}{j}(-1)^{r+k}k^d, 
$$
where $S(n,k)$  denotes Stirling number of the second kind$;$ see Section 5. 
\par 
In particular, 
$$
\ehk(A) + \mhk(A) = \dfrac{r! \cdot S(d,r)+s! \cdot S(d,s)}{d!}. 
$$
\endproclaim 
%
%
%
\head 
{\bf 1. Definition of the minimal Hilbert-Kunz multiplicity} 
\endhead 
\par 
In this section, we define the notion of 
{\it minimal Hilbert-Kunz multiplicity} and give its fundamental properties. 
In the following, let $A$ be a Noetherian excellent reduced local ring of 
positive characteristic $p > 0$ with perfect residue field $k = A/\frm$, 
unless specified. 
\medskip 
\subhead 
1A. Peskine-Szpiro functor  
\endsubhead 
\par 
First, let us recall the definition of Peskine--Szpiro functor. 
Let $\eA$ denote the ring $A$ viewed as an $A$-algebra 
via $F^e : A \to A \;(a \mapsto a^{p^e})$. 
Then $\FAe(-) = \eA \otimes_A -$ is a covariant functor from $A$-modules 
to $\eA$-modules. 
Since $\eA$ is isomorphic to $A$ as rings (via $F^e$), 
we can regard $\FA^e$ as a covariant functor from $A$-modules to themselves. 
We call this functor $\FAe$ {\it the Peskine--Szpiro functor\/} of $A$. 
The $A$-module structure on $\FAe(M)$ is such that $a'(a\otimes m) = a'a\otimes m$. 
On the other hand, $a'\otimes am = a'a^q \otimes m$. 
See e.g \cite{PS}, \cite{Ro} or \cite{Hu}. 
Suppose that an $A$-module $M$ has a finite 
presentation $A^m \overset \phi\to\longrightarrow A^n \to M \to 0$ where 
the map $\phi$ is defined by a matrix $(a_{ij})$.   
Then $\FAe(M)$ has a finite presentation 
$A^m \overset \phi_q\to\longrightarrow A^n \to \FAe(M) \to 0$ where 
the map $\phi_q$ is defined by the matrix $(a_{ij}^q)$.   
For example, $\FA^e(A/I) = A/I^{[p^e]}$, where 
$I^{[p^e]}$ is the ideal generated by $\{a^{p^e} : a \in I\}$. 
\subhead 
1B. Tight closure, Hilbert-Kunz multiplicity  
\endsubhead 
\par 
Using the Peskine--Szpiro functor, we define the notion of tight closure. 
\definition {Definition 1.1}  (\cite{HH1, HH2, Hu})
(i) Let $M$ be an $A$-module, and let 
$N$ be an $A$-submodule of $M$. 
Put $N^{[p^e]}_M = \Ker(\FA^e(M) \to \FA^e(M/N))$. 
Also, we denote by $F^e(x)$  
the image of $x$ the Frobenius map $M \to \FA^e(M) \; (x \mapsto 1\otimes x)$. 
Then  the {\it tight closure\/} $N^{*}_M$ of $N$ (in $M$)   
is the submodule generated by elements for which 
there exists an element $c \in A^0:= A \setminus 
\bigcup_{P \in \Min(A)} P$ such that for all sufficiently 
large $q = p^e$, $cF^e(x) \in N^{[q]}_M$. 
By definition, we put $I^{*} = I^{*}_A$. 
Also, we say that $N$ is {\it tightly closed \/} (in $M$) if $N^{*}_M =N$. 
\par 
(ii) A local ring $A$ in which every ideal is  
tightly closed is called {\it weakly F-regular}. 
Also, the ring whose localization is always weakly F-regular
is called {\it F-regular}.   
\par
(iv) A reduced ring $A$ is said to be {\it strongly F-regular\/} 
if for any element $c \in A^{0}$ 
there  exists $q =p^e$ such that the $A$-linear map 
$A \to A^{1/q}$ defined by $a \to c^{1/q}a$ is split injective. 
\par 
(v) A Noetherian ring $R$ is (resp. weakly, strongly) F-regular
if and only if so is $R_{\frm}$ for every maximal ideal $\frm$. 
\enddefinition 
\remark {Remark} 
Strongly F-regular rings are F-regular. 
In general, it is not known whether the converse is true or not. 
But it is known that F-finine $\bbQ$-Gorenstein weakly F-regular rings are always 
strongly F-regular; see \cite{AM, Mc, Wi}. 
\endremark
\medskip 
The notion of Hilbert-Kunz multiplicity plays the central role in this paper. 
\definition {Definition 1.2}{\rm (\cite{Ku2, Mo, Se1})}
Let $I$ be an $\frm$-primary ideal in $A$ and $M$ a finite $A$-module. 
Then we define the Hilbert-Kunz multiplicity $\ehk(I,\,M)$ of $M$ with 
respect to $I$ as 
$$
\ehk(I,\,M) := \lim_{e \to \infty} \frac{l_A(M/I^{[p^e]}M)}{p^{de}}. 
$$
By definition, we put $\ehk(I):=\ehk(I,A)$ and $\ehk(A) :=\ehk(\frm)$. 
\par 
Also, the multiplicity of $I$ is defined as 
$$
e(I) = \lim_{n \to \infty} \dfrac{d! \cdot l_A(A/I^n)}{n^d}. 
$$
\enddefinition 
\par
Suppose that $\widehat{A}$ (the $\frm$-adic completion of $A$) is reduced. 
Let $I \subseteq I'$ be $\frm$-primary ideals in $A$. 
Then it is known that $I'$ and $I$ have the same integral closure 
(i.e. $\overline{I'} = \overline{I}$) if and only if $e(I) = e(I')$. 
The similar result holds for tight closures and the Hilbert-Kunz multiplicities. 
\proclaim {Lemma 1.3} {\rm (cf. \cite{HH2, Theorem 8.17})}
Let $I \subseteq I'$ be $\frm$-primary ideals in $A$. 
\roster 
\item If $I' \subseteq I^{*}$, then $\ehk(I) = \ehk(I')$. 
\item 
Further assume that $\widehat{A}$ is excellent, equidimensional, and reduced. 
Then the converse of $(1)$ is also true. 
\endroster 
\endproclaim 
\subhead 
1C. Minimal Hilbert-Kunz multiplicity 
\endsubhead 
\par 
Our work is motivated by the following question. 
\proclaim {Question 1.4}
What is the minimal value of the difference 
$\ehk(I) - \ehk(I')$ for ideals $I' \supseteq I$ with $l_A(I'/I) =1$? 
\endproclaim 
In order to represent the \lq\lq difference'', we define the following notion.  
\definition {Definition 1.5 (Relative Hilbert-Kunz multiplicity)}
Let $L$ be an $A$-module,  
and let $N \subseteq M$ be  finite $A$-submodules 
of $L$ with $\Supp_A(M/N) \subseteq \{\frm\}$. 
Then we set  
$$
\ehk(N,\,M\,;L) := \liminf_{e \to \infty} 
\frac{l_A(M^{[p^e]}_L/N^{[p^e]}_L)}{p^{de}}. 
\tag"\bf (1.5.1)"
$$
We call $\ehk(N,\,M\,;L)$ the {\it relative Hilbert-Kunz multiplicity} with
respect to $N \subseteq M$ of $L$. 
In particular, 
$\ehk(I,\,I'\,;A) = \ehk(I) - \ehk(I')$ 
for $\frm$-primary ideals $I \subseteq I'$ in $A$. 
\par 
We further define the notion of the 
{\it relative minimal Hilbert-Kunz multiplicity} of $A$ 
as follows$:$ 
$$
\relmhk(A) := \inf\{\ehk(I,I';A)\,:\, I \subseteq I' \;\text{such that} \; 
l_A(I'/I) =1\}.   
\tag"\bf (1.5.2)"
$$
\enddefinition 

\medskip
The following properties
of the \lq\lq relative minimal Hilbert-Kunz multiplicity''
essentially follows from \cite{WY1,Theorem 1.5}. 

\proclaim {Proposition 1.6}
Let $A$ be an excellent local ring. Then  
\roster 
\item $0 \le \relmhk(A) \le 1$. 
\item $\relmhk(A) =1$ if and only if $A$ is regular. 
\item If $\relmhk(A) > 0$, then $A$ is weakly F-regular. 
\endroster 
\endproclaim 
\demo{\quad \bf Proof}
We first prove (3). 
Suppose that $A$ is not weakly F-regular. 
Then there exists an $\frm$-primary 
ideal $I$ such that $I \ne I^{*}$. 
Taking an ideal $I'$ such that 
$I \subseteq I' \subseteq I^{*}$ and $l_A(I'/I) =1$, 
we have $\ehk(I) = \ehk(I')$ by Lemma 1.3(1). 
This implies that $\relmhk(A) =0$; this is a contradiction. 
\par
To see (1),(2), it is enough to show that if $\relmhk(A) \ge 1$ then 
$A$ is regular. 
Actually, if $A$ is regular, 
then $\ehk(I) = l_A(A/I)$ for every ideal $I$ of $A$. 
Thus $\relmhk(A) =1$. 
\par 
Now suppose that $\relmhk(A) \ge 1$. 
Then $A$ is weakly F-regular, and thus is Cohen--Macaulay 
(cf. \cite{HH3}). 
Let $J$ be a parameter ideal of $A$. 
Then $\ehk(J) = e(J) = l_A(A/J)$. 
By the assumption that $\relmhk(A) \ge 1$, we get 
$$
\ehk(\frm) \le \ehk(J) - \ehk(\frm/J) = l_A(A/J) - l_A(\frm/J) =1. 
$$
Hence $A$ is regular by \cite{WY1,\,Theorem 1.5}. \oldqed 
\enddemo

\proclaim {Question 1.7}
Is the converse of Proposition $1.6(3)$ true?
\endproclaim 

\medskip
\remark {Remark}
In Section 3, we will give an affirmative answer to this question 
in case of $\bbQ$-Gorenstein F-regular local rings.  
\endremark 
\medskip 
In the following, 
to study $\relmhk(A)$ in detail, we introduce the notion of 
minimal Hilbert-Kunz multiplicity.  
\definition{Definition 1.8 (Minimal Hilbert-Kunz multiplicity)}
Let $E_A:=E_A(A/\frm)$ denote the injective hull of the residue field $k=A/\frm$. 
Let $z$ be a generator of the socle $\Soc (E_A) := \{x \in E_A\,|\, \frm x =0\}$ 
of $E_A$.  
Then we put 
$$
\mhk(A) :=\ehk(0,\Soc(E_A)\,; E_A)
= \liminf_{e \to \infty}\; \frac{l_A(A/\ann_A(F_A^e(z)))}{p^{ed}}, 
$$
where $F_A^e : E_A  \to \FAe(E_A)\;(u \mapsto 1\otimes u)$.
We call $\mhk(A)$ {\it the minimal Hilbert-Kunz multiplicity\/} of $A$. 
\enddefinition 
\proclaim {Example 1.9} 
Let $A = k[[x^e,\,x^{e-1}y,\,\ldots,xy^{e-1},\,y^e]]$, 
the $e$th Veronese subring of $k[[x,y]]$. 
Then $\relmhk(A)=\mhk(A) = 1/e$. 
\endproclaim 
\medskip
In the following proposition, we will prove 
$\relmhk(A)\ge \mhk(A)$. 
We expect that $\relmhk(A)=\mhk(A)$ always holds, but have no proof 
yet in general. 
This will be proved 
for Gorenstein local rings in the next section. 

\proclaim {Proposition 1.10} 
$\relmhk(A) \ge \mhk(A)$. 
\endproclaim 
\demo{\quad \bf Proof}
Since $\mhk(A) = \mhk(\widehat{A})$, we may assume that $A$ is complete. 
For given $\frm$-primary ideals $I \subseteq I'$ with $l_A(I'/I) = 1$,  
we must show that $\ehk(I) - \ehk(I') \ge \mhk(A)$. 
To see this, we need the following lemma. 
\proclaim {Lemma 1.11} 
Let $I \subseteq I'$ be $\frm$-primary ideals in $A$ such that $l_A(I'/I) =1$. 
Also, let $z$ be a generator of $\Soc(E_A)$. 
Then we can take elements $a \in I'$ and $u \in [0:_E I]$ which 
satisfy the following conditions$:$ 
\roster 
\item"(i)" $I' = I + aA$.  
\item"(ii)" $[0:_E I]= [0:_E I']+ Au$.  
\item"(iii)" $z = au$. 
\item"(iv)" $\frm a \subseteq I$. 
\item"(v)" $\frm u \subseteq [0:_E I']$. 
\endroster 
\endproclaim 
\demo{\quad \bf Proof} 
Since $l_A(I'/I) =1$, we can write $I'=I +aA$ for 
some $a \in I'\setminus I$ with $\frm a \subseteq I$. 
Applying Matlis duality to the natural surjective $A$-homomorphism 
$A/I \to A/I'$, we get $[0:_EI'] \subseteq [0:_EI]$
with $l_A([0:_EI]/[0:_EI']) = l_A(A/I) - l_A(A/I')=1$.  
Thus we can write $[0:_EI] = Au + [0:_EI']$ for some 
$u \in [0:_EI] \setminus [0:_EI']$ with $\frm u \subseteq [0:_EI']$. 
Then $Az = \Soc(E_A)$ is contained in $I'[0:_EI]$ 
since $I'[0:_EI] =a[0:_EI] \ne 0$. 
Hence we can write $z = rau$ for some $r \in A$. 
Then $r$ is a unit of $A$. 
Otherwise, since $ru \in [0:_EI']$, 
we have $z = a(ru) = 0$; this is a contradiction. 
Hence we may assume that $r=1$ without loss of generality. \oldqed   
\enddemo 
\par 
We now return to the proof of Proposition 1.10.  
Let $a \in I'$ and $u \in [0:_EI]$ be the elements described 
as in Lemma 1.11. 
We want to show the following claim. 
\demo{\bf Claim}
$I^{[q]}:a^q \subseteq \ann_A (F_A^e(z))\quad $ for 
all $q = p^e,\,e \ge 1$.   
\enddemo 
\par 
Let $c \in I^{[q]}:a^q$. 
Then $cF_A^e(z) = cF_A^e(au) = ca^q F_A^e(u)=0$, where 
the last equality follows from $I^{[q]}F_A^e([0:_EI])=0$.   
Hence $c \in \ann_A(F_A^e(z))$, as required. 
\medskip 
On the other hand, since $I^{'[q]} = I^{[q]} +Aa^q$, we have 
$$
I^{'[q]}/I^{[q]} = \frac{I^{[q]}+a^qA}{I^{[q]}} \cong A/I^{[q]}:a^q. 
$$
By the above claim, we get 
$$
\frac{l_A (I^{'[q]}/I^{[q]})}{q^d} = 
\frac{l_A (A/I^{[q]}:a^q)}{q^d}  
\ge \frac{l_A(A/\ann_A(F_A^e(z)))}{q^d}. 
$$ 
This yields the required inequality. \qed 
\enddemo 
\demo{Discussion {\rm 1.12}}
Let $I \subseteq I'$ be $\frm$-primary ideals in $A$ with $l_A(I'/I) =1$.  
Then an exact sequence $
\FA(A/\frm)  = A/\frm^{[q]}  \to A/I^{[q]} \to A/I^{'[q]} \to 0$ 
implies that 
$$
\ehk(I) - \ehk(I') \le \ehk(A).  
\tag"{\bf (1.12.1)}" 
$$
Thus 
$$
\max\{\ehk(I,\,I',A) \,:\, I \subseteq I' \;\text{with}\; l_A(I'/I) =1\} 
= \ehk(A). 
$$
\enddemo 
\remark {Remark}
One can prove the similar result as in Proposition 1.10 and Discussion 1.12
for any pair $(N,\,M)$ of $A$-submodules of an $A$-module 
$L$ with $N \subseteq M$ and $l_A(M/N) =1$.  
\endremark 
\medskip 
\head
{\bf 2. Gorenstein local rings}
\endhead 
In this section, we prove that  if $(A,\frm)$ is a Gorenstein local ring then 
$\ehk(J) - \ehk(J:\frm)$ is independent on the choice of parameter ideal $J$ of $A$. 
In fact, this invariant is equal to $\mhk(A)$ defined 
in the previous section. 
\medskip 
In the following, let $(A,\,\frm,\,k)$ be a $d$-dimensional 
local ring of characteristic $p >0$ and assume that 
$k$ is an infinite field.  
\par 
Now suppose that $A$ is Cohen-Macaulay. 
Then the highest local cohomology $H_{\frm}^d(A)$ may be identified with 
$\varinjlim A/(a_1^n,\ldots,a_d^n)A$, where  
$a_1,a_2,\ldots,a_d$ is a system of parameters for $A$ and the maps 
in the direct limit system are given 
by multiplication by $a = \prod_{i=1}^d a_i$. 
Put $(\underline{a})^{[n]} = (a_1^n,a_2^n,\ldots,a_d^n)A$ for 
every integer $n \ge 1$. 
Then any element $\eta \in H_{\frm}^d(A)$ can be represented 
as the equivalence class $[x + (\underline{a})^{[n]}]$ for 
some $x \in A$ and some integer $n \ge 1$. 
Also, note that we can write $(\underline{a})^{[q]} = J^{[q]}$ for all $q=p^e$. 
\medskip 
Considering the Frobenius action to $H_{\frm}^d(A)$, we can regard as 
$$
\FA^e(H_{\frm}^d(A)) \cong \varinjlim A/(\underline{a})^{[nq]} =H_{\frm}^d(A), 
$$
where $q =p^e$.  
Then $F^e(\eta) = [x^q \mod (\underline{a})^{[nq]}] \in H_{\frm}^d(A)$ for 
 $\eta = [x+(\underline{a})^{[n]}] \in H_{\frm}^d(A)$. 
See \cite{Sm} for more details.  
Using this fact, we can prove  the following theorem. 
\proclaim {Theorem 2.1}
Let $(A,\frm)$ be a Gorenstein local ring of characteristic $p > 0$. 
Then for any $\frm$-primary ideal 
$J$ of $A$ such that $\pd_A A/J < \infty$ and $A/J$ is Gorenstein, 
we have 
$$
\ehk(J) - \ehk(J:\frm) = \mhk(A). 
$$
In particular, $\relmhk(A) = \mhk(A)$. 
\endproclaim 
\demo{\quad \bf Proof}
First, we consider the case of parameter ideals. 
Since $A$ is Gorenstein, $E_A \cong H_{\frm}^d(A)$. 
In the above notation, the generator $z$ of $\Soc(E_A)$ can be written 
as $z = [b + J]$, where $b$ is a generator of 
$\Soc(A/J)$.  
For any element $c \in A$ and for all $q = p^e$, 
$$
cF_A^e(z)=cF_A^e([b+J])= [cb^q + J^{[q]}] = 0 \in H_{\frm}^d(A)
$$
if and only if there exists an integer $n \ge 1$ such that 
$$
cb^q \in (a_1^{nq},\ldots,a_d^{nq}) :  
(a_1^{n-1}\cdots a_d^{n-1})^q = J^{[q]}. 
$$
It follows that $\ann_A F^e_A(z) = J^{[q]}: b^q $. 
Hence we get 
$$
\mhk(A) 
 = \lim_{e \to \infty} \frac{l_A(A/J^{[q]}:b^q)}{q^d}
 = \lim_{e \to \infty} \frac{l_A((J:\frm)^{[q]}/J^{[q]})}{q^d}  
 =\ehk(J)-\ehk(J:\frm),   
$$
as required. 
\medskip
Next we consider the general case. 
Let $J$ be an $\frm$-primary ideal such that $\pd_A A/J < \infty$ and $A/J$ is 
Gorenstein. 
Take a parameter ideal $\frq$ which is contained in $J$. 
Then it is enough to show that 
$$
\ehk(J) - \ehk(J:\frm) = \ehk(\frq) -\ehk(\frq:\frm). 
$$
This follows from the argument as in the proof of \cite{Vr, Proposition 3.5}, 
but we put a sketch here for sake of completeness. 
As $\frq \subseteq J$, there exists a natural surjective map $A/\frq \to A/J$. 
Also, since both $A/\frq$ and $A/J$ are Gorenstein rings of dimension $d$, 
we have the following commutative diagram$:$ 
$$
\CD
 0 @. \to A=F_d'  @. \to \,\cdots\,\to @. A^d @>>>  A @>>> @. A/\frq @.\to 0 \\
    @.   @V \delta VV       @.           @VVV     @VVV @. @VV \text{nat} V 
@. @. \\
 0 @. \to A=F_d  @. \to \,\cdots\,\to @. A^n @>>>  A @>>> @. A/J @.\to 0 \\
\endCD
\tag"\bf (2.1.1)"
$$ 
where the horizontal sequences are minimal free resolutions of $A/\frq$ and  
$A/J$, respectively. 
In particular, the map $F_d' \to F_d$ is given by the multiplication of an 
element (say $\delta$). 
Then we have $J = \frq : \delta$. 
Let $L = J:\frm$ and $L'= \frq:\frm$. 
Since $\frq : (\frq,\,\delta L) = \frq : \delta L 
= (\frq:\delta) : L = J:L = \frm$, 
we get 
$$
(\frq,\delta L) = \frq : (\frq : (\frq,\,\delta L)) = \frq :\frm = L'. 
$$
Taking a Frobenius power $(\quad)^{[q]}$, 
we have $(\frq^{[q]},\,\delta^q L^{[q]}) = L^{'[q]}$. \par
On the other hand, $J^{[q]}  = \frq^{[q]} : \delta^q$ because if one can apply 
the Peskine--Szpiro functor to the diagram (2.1.1) then 
one can obtain the minimal 
free resolutions of $A/J^{[q]}$ and $A/\frq^{[q]}$, respectively. 
Hence 
$$
J^{[q]} : L^{[q]} = (\frq^{[q]} : \delta^q) : L^{[q]}
= \frq^{[q]} : (\frq^{[q]},\, \delta^qL^{[q]})
= \frq^{[q]} : L^{'[q]}. 
$$
If we write $L = J+aA$ and $L' = \frq+bA$, then 
$J^{[q]} : a^q  = \frq^{[q]} : b^q$. 
Thus 
$$
l_A(L^{[q]}/J^{[q]}) = l_A(A/J^{[q]}:a^q) 
= l_A(A/\frq^{[q]}:b^{q}) = l_A(L^{'[q]}/\frq^{[q]}) 
$$
for all $q = p^e$. 
The required assertion easily follows from this. 
\oldqed 
\enddemo 
\proclaim {Corollary 2.2} 
Under the same notation as in Theorem $2.1$, we have  
\roster 
\item $0 \le \mhk(A) \le 1$. 
\item $\mhk(A) = 1$ if and only if $A$ is regular. 
\item $\mhk(A) >0$ if and only if $A$ is F-regular. 
\endroster 
If, in addition, $e(A) \ge 2$, 
then 
$$
\mhk(A) \le \frac{e(A)-\ehk(A)}{e(A)-1}, 
$$ 
where $e(A) = e(\frm)$ denotes the usual multiplicity of $A$. 
\endproclaim 
\demo{\quad \bf Proof} 
The first two statement immediately follows from 
Proposition 1.6 and the previous theorem. 
\par 
To see (3), it is enouth to see \lq\lq if'' part.  
Actually, since Gorenstein weakly F-regular local ring is always F-regular, 
\lq\lq only if'' part follows from Proposition 1.6 and Theorem 2.1. 
Now suppose that $A$ is weakly F-regular, that is, 
every ideal of $A$ is tightly closed. 
Thus $\ehk(J) \ne \ehk(J:\frm)$ for every parameter ideal $J$ of $A$ 
by Lemma 1.3(2). 
Therefore we have $\mhk(A) = \ehk(J) - \ehk(J:\frm) > 0$. 
\par 
To see the last inequality,  
taking a minimal reduction $J$ of $\frm$, we have 
$$
\ehk(J) - \ehk(\frm) \ge l_A(\frm/J) \cdot \mhk(A).  
$$
This yields a required inequality since $\ehk(J) = e(J) = e(A)$. 
\oldqed 
\enddemo 
\remark {Remark}
In \cite{HuL}, Huneke and Leuschke independently proved the similar result. 
In fact, in case of Gorenstein local rings, the notion of minimal 
Hillbert-Kunz multiplicity conincides with that of 
\lq\lq {\it rational signature}'' which was defined in \cite{HuL}. 
Also, see \cite{AL} for more details. 
\endremark 

\proclaim {Example 2.3}
Assume that $A$ is a hypersurface local ring of multiplicity $2$. 
Then we have $\mhk(A) = 2- \ehk(A)$. 
\endproclaim 
\medskip 
Let $A$ be a two-dimensional Gorenstein F-regular local ring 
which is not regular. 
Then $e(A) =2$ since $A$ has minimal multiplicity.
Moreover, suppose that $k$ is an algebraically closed field. 
Then it is known that the $\frm$-adic completion $\widehat{A}$ of $A$ is 
isomorphic to the completion of the invariant subring by a finite subgroup $G 
\subseteq SL(2,k)$ which acts on the polynomial ring $k[x,y]$.  
Furthermore, we have $\ehk(A) = 2 - 1/|G|$; see \cite{WY1,Theorem 5.1}. 
Hence $\mhk(A) = 1/|G|$ by Example 2.3. 
This result will be generalized in Section 4. 

By the above observation, we have an inequality 
$\mhk(A) \le \frac{1}{2}$ for hypersurface local rings with 
$\dim A = e(A) =2$. 
We can extend this result for hypersurface local rings of 
higher dimension. 

\proclaim {Proposition 2.4}
Let $(A,\frm)$ be a hypersurface local ring of chacteristic $p>0$. 
Suppose that $e(A) = \dim A = d \ge 1$. 
Then 
$$
\mhk(A) \le \dfrac{1}{2^{d-1}\cdot (d-1)!}.
$$
\endproclaim

\demo{\quad \bf Proof}
We may assume that $A$ is a complete F-rational local domain with infinite 
residue field. 
Let $J$ be a minimal reduction of $\frm$. 
Take an element $x \in \frm$ such that $\frm = xA +J$. 
Then since $x^{d-1}$ is a generator of $\Soc (A/J)$, we have 
$$
\mhk(A) = \lim_{e \to \infty} \frac{l_A(Ax^{(d-1)q} + J^{[q]}/J^{[q]})}{q^d}
$$
by Theorem 2.1. 
For any $q=p^e$, we have the following claim. 
\demo{\bf Claim}
$$
l_A(Ax^{(d-1)q}+J^{[q]}/J^{[q]}) \le 
2 \cdot l_A(A/\frm^{\lfloor \frac{q+1}{2}\rfloor}).
$$
\enddemo
To prove the claim, we put $B = A/J^{[q]}$, $y = x^{(d-1)q}$ and 
$\fra = \frm^{\lfloor \frac{q+1}{2}\rfloor}$. 
Then since $y\fra^2 \subseteq x^{(d-1)q}\frm^q \subseteq \frm^{dq} \subseteq J^{[q]}$, 
we have $y\fra B \subseteq 0 : \fra B=K_{A/\fra}$. 
By Matlis duality, we get 
$$
l_A(y B) \le l_A(yB/y\fra B) + l_A(y\fra B) \le l_A(A/\fra) + l_A(A/\fra), 
$$
as required. 
Since $l_A(A/\frm^n) = \frac{e(A)}{d!} n^d + O(n^{d-1})$ for all 
large 
enough $n$, the assertion easily follows from the claim. 
\oldqed
\enddemo

\demo{Discussion 2.5}
Let $(A,\frm)$ be a three-dimensional F-regular hypersurface local ring. 
Then $\ehk(A) \ge \frac{2}{3}e(A)$ by the following formula$:$ 
$$
\ehk(A) \ge \frac{e(A)}{\pi} \int_{-\infty}^{\infty} 
\left(\frac{\sin \theta}{\theta}\right)^{d+1} \,d\theta
\;= \;\frac{e(A)}{2^d d!}
\sum_{i=0}^{[\frac{d}{2}]}(-1)^i (d+1-2i)^d {d+1 \choose i}.  
$$
In particular, if futhermore $e(A) =3$, then $\ehk(A) \ge \frac{2}{3}\cdot 3 = 2$. 
Thus $\mhk(A) \le \frac{3-2}{3-1} = \frac{1}{2}$ by Corollary 2.2.
On the other hand, Proposition 2.4 implies that $\mhk(A) \le \frac{1}{8}$. 
\enddemo 

\proclaim {Question 2.6}
Let $d$ be an integer with $d \ge 2$, and 
let 
$$
A = k[[x_0,x_1,\ldots,x_d]]/(x_0^d+x_1^d+\cdots + x_d^d),
$$
where $k$ is a field of charcteristic $p>0$. 
Does $\mhk(A) = \frac{1}{2^{d-1}(d-1)!}$ hold if $p > d$?
\endproclaim



\medskip 
\head 
{\bf 3. Canonical covers}
\endhead 
\par
In the previous section, we have showed how to compute $\mhk(A)$ 
in the case of Gorenstein local rings. 
In this section, we study the minimal Hilbert-Kunz multiplicity in case 
of $\bbQ$-Gorenstein local rings using \lq\lq canonical cover''. 
\par
Now let us recall the notion of canonical cover. 
Let $A$ be a normal local ring 
and let $I$ be a divisorial ideal 
(i.e. an ideal of pure height one) of $A$.  
Also, let $\Cl(A)$ denote the divisor class group of $A$.
Now suppose that $\cl(I)$ is a torsion element in $\Cl(A)$, that is, 
$I^{(r)} := \cap_{P \in \Ass_A(A/I)} I^rA_P \cap A$ is a principal ideal 
for some integer $r \ge 1$. 
Putting $r = \ord(cl(I))$, one can write as $I^{(r)} = fA$ for 
some element $f \in A$. 
Then a $\bbZ^r$-graded $A$-algebra
$$
B(I,r,f):= A \oplus I \oplus I^{(2)} \oplus \cdots \oplus I^{(r-1)}
= \sum_{i=0}^{r-1} I^{(i)}t^i, \quad \text{where}\;\; t^rf=1 
$$
is called the $r$-{\it cyclic cover} of $A$ with respect to $I$. 
Also, suppose that $r$ is relatively prime to $p=\chara(A) >0$.  
Then $B(I,r,f)$ is a local ring with the unique maximal ideal 
$\frn := \frm \oplus I \oplus \cdots \oplus I^{(r-1)}$, and 
the natural inclusion $A \hookrightarrow  B(I,r,f)$ is \'etale 
in codimension one; thus $B(I,r,f)$ is also normal. 
\par 
Also, we further assume that $A$ admits a canonical module $K_A$. 
Note that one can regard $K_A$ as an ideal of pure height one. 
The ring $A$ is called $\bbQ$-Gorenstein if 
$\cl(K_A)$ is a torsion element in $\Cl(A)$. 
Put $r :=\ord(\cl(K_A)) < \infty$. 
Then the $r$-cyclic cover with respect to $K_A$   
$$
B :=  A \oplus K_A \oplus K_A^{(2)} \oplus \cdots \oplus K_A^{(r-1)}
$$ 
is called the {\it canonical cover} of $A$. 
\par 
Roughly speaking, the notion of $\bbQ$-Gorenstein F-regular local rings 
(resp. Gorenstein F-rational local rings)
is ``a positive characateristic analogy'' 
of that of {\it log terminal singularities} 
(resp. $1$-Gorenstein rational singularity) in characteristic zero. 
Thus \lq\lq canonical cover trick'' is one of the important tools 
not only in the theory of singularities 
but also in the theory of tight closures. 
Actually, if $A$ is a $\bbQ$-Gorenstein F-regular local ring then 
the canonical cover is not only F-regular but also Gorenstein. 
Further, $A$ is a direct summand of $B$ as an $A$-module. 
Thus one can reduce several problems of $\bbQ$-Gorenstein 
F-regular local rings to those of Gorenstein F-regular local rings
(this trick is called ``canonical cover trick''). 
See e.g. \cite{Wa3}, \cite{TW} for details. 
\par 
In fact, in our context, we can prove the following theorem. 

\proclaim {Theorem 3.1}
Let $(A,\,\frm,k)$ be a $\bbQ$-Gorenstein F-regular local ring and 
let $B =\oplus_{i=0}^{r-1}K_A^{(i)}t^i$ be the canonical cover of $A$,  
where $r$ is the order of $\cl(K_A)$ in $\Cl(A)$ and $t^rf=1$. 
Also, suppose that $r$ is relatively prime to $p = \chara(A)$.  
Then we have 
$$
\mhk(B) = r \cdot \mhk(A). 
$$
\endproclaim 
\medskip
The following corollary gives a partial answer to Question 1.7. 

\proclaim {Corollary 3.2}
Let $A$ be a $\bbQ$-Gorenstein F-regular local 
ring of characteristic $p >0$ 
such that $(\ord(\cl(K_A)),\,p)  =1$. 
Then $\mhk(A) > 0$. 
\endproclaim 

\medskip
In the following, we shall prove Thoorem 3.1. 
We begin with recalling several properties of  
canonical covers 
for convenience of the readers. 
\proclaim {Lemma 3.3} 
Let $(A,\,\frm,\,k)$ be a Cohen-Macaulay normal local ring, and 
suppose that $A$ is $\bbQ$-Gorenstein. 
Let $B = \oplus_{i=0}^{r-1}K_A^{(i)}t^i$ be the canonical cover of $A$, 
where $t^rf =1_A$. 
Then the following statements hold. 
\roster 
\item $B$ is quasi-Gorenstein, that is, $B \cong K_B$ as $B$-modules. 
In particular, $B$ is Gorenstein if it is Cohen-Macaulay. 
\item $A$ is strongly F-regular if and only if so is $B$. 
\endroster 
\par 
In the following, we further assume that $B$ is Cohen-Macaulay. 
\roster 
\item"(3)" The injective hull $E_B := E_B(B/\frn)$ of $B/\frn$ is given as follows$:$
$$
E_B = \bigoplus_{i=0}^{r-1} H_{\frm}^d(K_A^{(i)})t^i. 
$$
\item"(4)" $\Soc_B(E_B)= \Hom_B(B/\frn,E_B)$ is generated by $zt$,  
where $z$ is a generator of the socle of $E_A \cong H_{\frm}^d(K_A)$. 
\endroster 
\endproclaim 
\demo{\bf Proof}
The assertion (1) follows from \cite{TW, Sect.3}. 
Also, the assertion (2) follows from \cite{Wa3,Theorem 2.7}. 
\par 
In the following, assume that $B$ is Cohen-Macaulay. 
Then since $B$ is Gorenstein by (1) and $\frm B$ is $\frn$-primary, 
we have 
$$
E_B \cong H_{\frn}^d(B) \cong H_{\frm}^d(B) \cong
\bigoplus_{i=0}^{r-1} H_{\frm}^d(K_A^{(i)})t^i. 
$$
Thus we get the assertion (3). 
To see (4), it ie enough to show that $zt \in \Soc_B(E_B)$ 
since $\dim_k \Soc_B(E_B) =1$. 
Namely, we must show that $az = 0$ in $H_{\frm}^d(K_A^{(i+1)})$ 
 for all $i$ with $1 \le i \le r-1$ and for all $a \in K_A^{(i)}$. 
\par
Fix an integer $i$ with $1 \le i \le r-1$ 
and suppose that $0 \ne a \in K_A^{(i)}$.  
Applying the local cohomology functor to the short exact sequence 
$$
 0 \to K_A \overset a \to\longrightarrow 
K_A^{(i+1)} \to K_A^{(i+1)}/aK_A \to 0 
$$
implies that 
$$
0 = H_{\frm}^{d-1}(K_A^{(i+1)}) \to H_{\frm}^{d-1}(K_A^{(i+1)}/aK_A) 
\to H_{\frm}^d(K_A) \overset a \to\longrightarrow H_{\frm}^d(K_A^{(i+1)}), 
$$
where the first vanishing follows from the fact that $K_A^{(i+1)}$ is a 
direct summand of a maximal Cohen-Macaulay $A$-module $B$. 
To get the lemma, it is enough to show the following claim$:$ 
\demo{\bf Claim}
$H_{\frm}^{d-1}(K_A^{(i+1)}/aK_A) \ne 0$.
\enddemo 
\par 
Since $A$ is Cohen-Macaulay, $aK_A \cong K_A$ is 
a maximal Cohen-Macaulay $A$-module. 
Hence $aK_A$ is a divisorial ideal of $A$. 
Thus it is enough to show that $K_A^{(i+1)}/aK_A \ne 0$. 
Suppose not. 
Then we have $(i+1) \divi(K_A) = \divi(K_A) + \divi(a)$, and 
thus $i \cdot \cl(K_A) =0$. 
This contradicts the assumption that $r = \ord(\cl(K_A))$. 
Hence we get the claim, as required. \oldqed 
\enddemo 
\demo{\bf Proof of Theorem 3.1}
Let $x_1,\,\,\ldots,x_d$ be a system of parameters of $A$ and fix it. 
Since $A$ is Cohen--Macaulay, we have $E_A = H_{\frm}^d(K_A) 
= \varinjlim \, K_A/\underline{x}^{[q]}K_A$. 
Also, $A$ is $\bbQ$-Gorenstein, 
one can regard the Frobenius map $\FAe$ in $E_A$ as 
$$
F_A^e : E_A \to \FAe(E_A) 
\cong H_{\frm}^d(K_A^{(q)}) = \varinjlim K_A^{(q)}/\underline{x}^{[n]}K_A^{(q)}
\; \left([b+\underline{x}K_A] \to 
[b^q+\underline{x}^{[q]}K_A^{(q)}]\right); 
$$
see \cite{Wa3} for details. 
Thus we have 
$$
\mhk(A) = \liminf_{q \to \infty} \;
l_A\left(\dfrac{z^qA + \underline{x}^{[q]}K_A^{(q)}}
{\underline{x}^{[q]}K_A^{(q)}}\right)
\bigg/ q^d, 
\tag"\bf (3.1.1)"
$$
where we denote the inverse image of the generator $z$ of the socle 
of $K_A/\underline{x}K_A$ by the same symbol as $z$.  
\par 
On the other hand, since $zt \in K_At$ generates the socle of $E_B$ by Lemma 3.3, 
we get  
$$
\mhk(B) = \lim_{q \to \infty} 
l_A\left(\dfrac{z^qt^q B + \underline{x}^{[q]}B}{\underline{x}^{[q]}B}\right)
\bigg/q^d
\tag"\bf (3.1.2)"
$$
by Theorem 2.1. 
Also, as $B$ is a $\bbZ/r\bbZ$-graded ring 
(especially, $K_A^{(i+r)}t^{i+r} = K_A^{(i)}t^i$), (3.1.2) is reformulated 
as follows$:$
$$
\mhk(B) = \sum_{i=0}^{r-1} 
\lim_{q \to \infty} 
l_A\left(\dfrac{z^qK_A^{(i)} + \underline{x}^{[q]}K_A^{(i+q)}}
{\underline{x}^{[q]}K_A^{(i+q)}}\right)\bigg/q^d
\tag"\bf (3.1.3)"
$$
If necessary, we may assume that $q \equiv 1 \pmod r$. 
Taking a nonzero element $a_i \in K_A^{(i)}$ for each $i$ 
with $0 \le i \le r-1$, 
we consider the following commutative diagram with exact rows$:$
$$
\CD
 0 \to @. K_q  @>>>
\dfrac{z^qA + \underline{x}^{[q]}K_A^{(q)}}{ \underline{x}^{[q]}K_A^{(q)}} 
@> a_i >> 
\dfrac{z^qK_A^{(i)} + \underline{x}^{[q]}K_A^{(i+q)}}{ \underline{x}^{[q]}K_A^{(i+q)}} 
@>>>
C_q 
@.\to 0  \\ 
@. @VVV  @VV inj. V  @VV inj. V  @VVV \\
 0 \to @. X_q  @>>>
\dfrac{K_A^{(q)}}{ \underline{x}^{[q]}K_A^{(q)}} 
@> a_i >> 
\dfrac{K_A^{(i+q)}}{ \underline{x}^{[q]}K_A^{(i+q)}} 
@>>>
Y_q 
@.\to 0.   
\endCD
$$ 
\par 
In order to complete the proof of the theorem, 
it suffices to prove the following claim. 

\demo{\bf Claim}
$\displaystyle{\lim_{q \to \infty}} \dfrac{l_A(K_q)}{q^d} = 
\displaystyle{\lim_{q \to \infty}} \dfrac{l_A(C_q)}{q^d} = 0$. 
\enddemo
\par 
First, note that if $N$ is a finite $A$-module with $\dim N \le d-1$ 
then $l_A(N/\underline{x}^{[q]}N)/q^d = 0$. 
By definition of $Y_q$, we have 
$Y_q = K_A^{(i+q)}/(a_iK_A^{(q)}+\underline{x}^{[q]}K_A^{(i+q)}) 
\cong \left(K_A^{(i+1)}/a_iK_A \right) \otimes_A A/\underline{x}^{[q]}$. 
Since $\dim K_A^{(i+1)}/a_iK_A \le d-1$, we get 
$\displaystyle{\lim_{q \to \infty}}l_A(Y_q)/q^d =0$. 
On the other hand, as $q \equiv 1 \pmod r$, we have 
$$
\aligned 
\lim_{q \to \infty} \dfrac{l_A(K_A^{(q)}/\underline{x}^{[q]}K_A^{(q)})}{q^d}
& = \ehk(\underline{x}) \cdot \rank_A K_A = \ehk(\underline{x}), \\
\lim_{q \to \infty} \dfrac{l_A(K_A^{(i+q)}/\underline{x}^{[q]}K_A^{(i+q)})}{q^d} 
& = \ehk(\underline{x}) \cdot \rank_A K_A^{(i+1)} = \ehk(\underline{x}). 
\endaligned
$$
That is, 
$\displaystyle{\lim_{q \to \infty}} l_A(X_q)/q^d = 
\displaystyle{\lim_{q \to \infty}} l_A(Y_q)/q^d = 0$ and 
thus $\displaystyle{\lim_{q \to \infty}} l_A(K_q)/q^d =0$. 
\par
On the other hand, 
$$
\aligned
C_q 
& = \dfrac{z^qK_A^{(i)}+\underline{x}^{[q]}K_A^{(i+q)}}
{a_iz^qA+\underline{x}^{[q]}K_A^{(i+q)}} 
\cong \dfrac{z^qK_A^{(i)}}{a_iz^qA+z^qK_A^{(i)} \cap \underline{x}^{[q]}K_A^{(i+q)}} \\
& = \dfrac{z^qK_A^{(i)}}{a_iz^qA+z^q [K_A^{(i)} 
\cap (\underline{x}^{[q]}K_A^{(i+q)}:z^q)]} \\
& \cong \dfrac{K_A^{(i)}}{a_iA+ [K_A^{(i)} \cap (\underline{x}^{[q]}K_A^{(i+q)}:z^q)]}.  
\endaligned
$$
Since $\frm^{[q]}K_A^{(i)} 
\subseteq K_A^{(i)} \cap (\underline{x}^{[q]}K_A^{(i+q)}:z^q)$ 
by the choice of $z \in K_A$, we get 
$$
l_A(C_q) \le l_A(K_A^{(i)}/a_iA+\frm^{[q]}K_A^{(i)})
= l_A(K_A^{(i)}/a_iA \otimes_A A/\frm^{[q]}). 
$$
By the similar argument as above, we can prove $\displaystyle{\lim_{e\to \infty}}
l(C_q)/q^d =0$, as required. 
\oldqed
\enddemo

\example {Question 3.4}
Let $A$ be a weakly F-regular local ring and let $I$ be a divisorial ideal of $A$
such that $\cl(I)$ has a finite order (say $r$). 
If $B = A \oplus It \oplus I^{(2)}t^2 \oplus \cdots \oplus I^{(r-1)}t^{r-1}$, 
the $r$-cyclic cover, then does $\mhk(B) = r \cdot \mhk(A)$ hold ?  
\endexample 

\example {Question 3.5}
Under the same notation as in Theorem $3.1$, 
does $\relmhk(B) = r \cdot \relmhk(A)$ hold? 
Equivalently, does $\relmhk(A) = \mhk(A)$ hold?
\endexample 

\example {Discussion 3.6}
Let $(A,\frm) \hookrightarrow (B,\frn)$ be a module-finite extension of local domains 
of characteristic $p >0$. 
Put $r = [Q(B):Q(A)]$. 
If $A \hookrightarrow B$ is pure, then $\relmhk(B) \le r \cdot \relmhk(A)$. 
\par
In fact, let $I\subseteq I'$ be $\frm$-primary ideals of $A$ with 
$l_A(I'/I) =1$. 
Since $A \subseteq B$ is pure, we have $L = LB \cap A$ for 
all $L \subseteq A$. 
In particular, $IB \ne I'B$. 
By definition of $\relmhk(B)$, we get
$$
\relmhk(B) \le \ehk(IB) - \ehk(I'B) = r \cdot \{\ehk(I) - \ehk(I')\}. 
$$
Hence $\relmhk(B) \le  r \cdot \relmhk(A)$.
\endexample

\medskip 
\head
{\bf 4. Quotient singularities} 
\endhead 
\par 
In this section, as an application of Theorem 3.1, 
we study the minimal Hilbert-Kunz multiplicities  for
quotient singularities (i.e. the invariant subrings by a finite group; see below 
in detail). 
In general, quotient singularities are not necessalily Gorenstein but 
$\bbQ$-Gorenstein normal domains. 
So using the notion of canonical cover trick, we can reduce our problem 
to the case of Gorenstein rings. 
\par 
Let $k$ be a field 
and $V$ a $k$-vector space of finite dimension (say $d = \dim_k V$). 
Assume that a finite subgroup $G$ of $GL(V) \cong GL(d,\,k)$ acts linearly 
on $S  := \Sym_k(V) \cong k[x_1,\ldots,x_d]$, a polynomial ring with 
$d$-variables over $k$. 
Then  
$$
S^G := \{f \in S \,:\, g(f) = f \quad \text{for all} \;g \in G\} 
$$
is said to be an {\it invariant subring} of $S$ by $G$.  
\par 
In this section, 
we consider positive characteristic (say $p=\chara(k)$) only,  
and assume that the order $|G|$ is non-zero in $k$, 
that is, $|G|$ is not divided by $p$. 
Then the existence of the Reynolds operator 
$$
\rho : S \to S^G \quad \left(\, a \mapsto |G|^{-1} \sum_{g \in G} g(a) \right),  
$$
claims that $S^G$ is a direct summand of $S$.  
Put $\frn = (x_1,\,\ldots,x_d)S$ and $\frm = \frn \cap S^G$. 
Then the ring $A = (S^G)_{\frm}$ is said 
to be a {\it quotient singularity} (by a finite group $G$). 
Such a quotient singularity as above is a $\bbQ$-Gorenstein strongly 
F-regular domain, but it is not always Gorenstein; 
see e.g. \cite{Wa1,\,Wa2} for details.  
\medskip
In \cite{WY1}, 
we gave a formula for Hilbert-Kunz multiplicity 
$\ehk(A)$ of quotient singularities as follows. 
\proclaim {Theorem 4.1} {\rm (cf.~\cite{WY1,~Theorem 2.7},~\cite{BCP})} 
Under the same notation as above, we have 
$$
\ehk(I) = \frac{~1~}{|G|}\, l_A(S_{\frn}/IS_{\frn}),  
$$
for every $\frm$-primary ideal $I$ in $A$. 
In particular, $\ehk(A) = \frac{~1~}{|G|}\, \mu_A(S_{\frn})$, 
where $\mu_A(M)$ denotes the number of minimal system of 
generators of a finite $A$-module $M$. 
\endproclaim 
\par 
The main purpose of this section is to prove the following theorem. 
\proclaim {Theorem 4.2} 
Let $A=(S^G)_{\frm}$ be a quotient singularity 
by a finite group $G$ described as above. 
Also, assume that $G$ contains no pseudo-reflections. 
Then we have 
$$
\mhk(A) = \frac{~1~}{|G|}.  
$$
\endproclaim 
\demo{\quad \bf Proof}
First, suppose that $G \subseteq SL(d,k)$. 
Then $S^G$ is Gorenstein by \cite{Wa1, Theorem 1a}. 
Since $G$ acts linearly on $S$, $S^G$ is a graded subring of $S$. 
Thus one can take a homogeneous system of parameters  
$a_1,\,\ldots,a_d$ of $S^G$ with the same degree $m$. 
Also, we may assume that $m$ is a multiple of $|G|$. 
Put $J=(a_1,\,\ldots,a_d)S^G$. 
Then since $S/JS$ is a homogeneous Artinian Gorenstein ring
having the same Hilbert function as that of $S/(x_1^m,\,\ldots,x_d^m)S$, 
there exists an element $z \in S_{d(m-1)}$ which generates $\Soc(S/JS)$. 
Then we have $z \in S^G$, which follows from the proof of
\cite{Wa1,Theorem 1a}. 
But since this fact is essential point in the proof in this case, 
we put a sketch here. 
\par
To see $z \in S^G$, it is enough to show that $z \in S^{\langle g \rangle}$ for any 
element $g \in G$. 
The property $z \in S^{\langle g \rangle}$ does not change if we consider 
$S \otimes_k \overline{k}$ instead of $S$, where $\overline{k}$ is the 
algebraic closure of $k$. 
Then we may assume $k=\overline{k}$ and assume $g$ is diagonal. 
Then $x_1\cdots x_d \in S^{\langle g \rangle}$ and 
$x_i^m \in S^{\langle g \rangle}$ since $\det(g)=1$. 
If we put $(\underline{x})^{[m]} = (x_1^m, \cdots, x_d^m)$, then 
$$
\dim_k [S^{\langle g \rangle}/JS^{\langle g \rangle}]_{d(m-1)} 
= \dim_k [S^{\langle g \rangle}/(\underline{x})^{[m]}
S^{\langle g \rangle}]_{d(m-1)} \ge 1.
$$
On the other hand, since $JS^{\langle g \rangle} 
= JS \cap S^{\langle g \rangle}$, we have 
$$
\dim_k [S^{\langle g \rangle}/JS^{\langle g \rangle}]_{d(m-1)} 
\le \dim_k [S/JS]_{d(m-1)} = 1. 
$$
It follows that $z \in S^G$, as required. 
\par 
Now let $J$, $z$ be as above. 
Then $JA:\frm A = (J,z)A$ and $JS: \frn = (J,z)S$. 
Hence we get 
$$
\align
\ehk(JA) - \ehk(JA:\frm A) 
& = \frac{~1~}{|G|} l_A(S_{\frn}/JS_{\frn}) - 
\frac{~1~}{|G|} l_A(S_{\frn}/(J:\frm)S_{\frn}) \cr 
& = \frac{~1~}{|G|} l_{S_\frn}(JS_{\frn}:\frn/JS_{\frn}) =
\frac{~1~}{|G|}. 
\endalign
$$
The required assertion follows from Theorem 2.1. 
\medskip 
Next, we consider the general case. 
If we put $H = G \cap SL(n,k)$, then 
$S^H$ is Gorenstein by \cite{Wa2, Theorem 1}. 
Further, since $H$ is a normal subgroup of $G$ and 
$G/H$ is a finite subgroup  of $k^{\times}$,  
$G/H$ is a cyclic group.  
Say $G/H = \langle \sigma H\rangle$ and $r = |G/H|$.  
Also, $S^G = (S^H)^{\langle \sigma \rangle}$. 
\par 
Then $B = (S^H)_{\frn \cap S^H}$ 
is a cyclic $r$-cover of $A=(S^G)_{\frm}$. 
In fact, it is known that $B$ is isomorphic to the canonical cover of $A$. 
$$
B \cong A \oplus K_A t \oplus K_A^{(2)}t^2 \oplus \cdots \oplus K_A^{(r-1)}t^{r-1},
$$
where $K_A^{(r)}=fA,\;t^rf =1$. See \cite{TW} in detail. 
\par
Since $\mhk(B) = \frac{~1~}{|H|}$, by Theorem 3.1, 
we get
$$
\mhk(A) = \frac{~1~}{r} \mhk(B) = \frac{~1~}{(G:H)|H|} = \frac{~1~}{|G|},  
$$
as required. \oldqed 
\enddemo 

\example {Conjecture 4.3} 
Under the same notation as in Theorem 4.2, 
$\relmhk(A) = 1/|G|$. 
\endexample 
\medskip 
\head 
{\bf 5. Segre products} 
\endhead 
\par 
Throughout this section, let $k$ be a perfect field of characteristic $p>0$, 
and let $R=k[x_1,\,\ldots,x_r]$ (resp. $S=k[y_1,\,\ldots,y_s]$) be a 
polynomial ring with $r$-variables (resp. $s$-variables) over $k$. 
Also, we regard these rings as homogeneous $k$-algebras with 
$\deg(x_i) = \deg(y_j)=1$ as usual. 
Then we define the graded subring $A = R \# S$ of $R \otimes_k S$ 
by putting $A_n := R_n \otimes S_n$ for all integer $n \ge 0$. 
Then $A = R \# S$ is said to be the {\it Segre product} of $R$ and $S$. 
Actually, the ring $A$ is the coordinate ring of the Segre 
Embedding $\bbP^{r-1} \times \bbP^{s-1} \hookrightarrow \bbP^{rs-1}$. 
\par 
As the Segre product $A$ is a direct summand of $R \otimes_k S$ 
(which is isomorphic to a polynomial ring with $r+s$-variables), 
it is a strongly F-regular domain. 
Further, it is known that $\dim A = r+s-1$ and $e(A) = {r+s-2 \choose r-1}$. 
See also \cite{GW, Chapter 4} for more details. 
\medskip 
The main purpose of this section is to compute the minimal Hilbert-Kunz 
multiplicity $\mhk(A)$ for Segre products. 
Before stating our result, we recall related results. 
\par 
In \cite{BCP}, Buchweitz, Chen and Purdue have given the Hilbert-Kunz multiplicity 
$\ehk(A)$ of $A$. 
Also, Eto and the second-named author \cite{Et} simplified their result
in terms of \lq\lq Stirling numbers of the second kind'' as follows. 

\proclaim {Theorem 5.1} 
{\rm (cf. \cite{BCP, 2.2.3}, \cite{EY, Theorem 3.3}, \cite{Et})} 
Suppose that $2 \le r \le s$ and put $d = r+s-1$. 
Let $A=k[x_1,\,\ldots,x_r]\#k[y_1,\,\ldots,y_s]$. 
Then 
$$ 
\ehk(A)  =\frac{s!}{d!}\,S(d,s) 
-\frac{1}{d!}\sum_{k=1}^{r-1} \sum_{j=1}^{r-k}
\binom{r}{k+j} \binom{s}{j}(-1)^{r+k}k^d, 
$$
where $S(n,k)$  denotes the Stirling number of the second kind$;$ see below. 
\endproclaim  

\medskip
Stirling numbers of the second kind also plays an important role 
in the study of the minimal Hilbert-Kunz multiplicity of the Segre product. 
So we recall the notion of Stirling numbers. 

\definition {Definition 5.2} 
{\rm (\cite{St,Chapter 1, \S 1.4})}
We denote by $S(n,k)$ the number of partitions of the 
set $[n]:=\{1,\,\ldots,n\}$ into $k$ blocks. 
Then $S(n,k)$ is called the {\it Stirling number of the second kind}. 
\enddefinition 

\medskip 
The following properties are well-known. See \cite{St}. 

\proclaim {Fact 5.3}
If we denote by $S(n,k)$ the Stirling number of the second kind, then
\roster 
\item The following identy holds$:$
$$
x^n=  \sum_{k=0}^{\infty} S(n,k) \;x(x-1)\cdots (x-k+1). 
$$
\item $S(n,k) = k S(n-1,k) + S(n-1,k-1)$ and $S(0,0)=1$. 
\item $S(n,k)$ admits the following exponential generating 
function$:$
$$
\sum_{n \ge k} S(n,k) \dfrac{x^n}{n!}  = \dfrac{1}{k!} (e^x-1)^k. 
$$
In particular, 
$$
S(n,k)  = \dfrac{1}{k!} \sum_{i=0}^k (-1)^{k-i} \binom{k}{i} i^n. 
$$
\endroster 
\endproclaim 
\medskip
For example, since $S(s+1,s) = \binom{s+1}{2}$, we have 
the following example. 

\example {Example 5.4}
Let $A = R \# S = k[x_1,\,x_2]\# k[y_1,\ldots,y_s]$, which is isomorphic 
to the Rees algebra $S[\frn t]$ over $S$. 
Then 
$$
\ehk(A) = s \left(\frac{1}{~2~}+\frac{1}{(s+1)!}\right). 
$$  
\endexample

\medskip 
In the following, we will give a formula for the minimal Hilbert-Kunz 
multiplicities of the Segre product. 
Now let $A$ be the Segre product of $R$ and $S$ described as above$:$ 
$A = R\#S =k[x_1,\,\ldots,x_r]\#k[y_1,\,\ldots,y_s]$, 
and suppose that $2 \le r \le s$. 
Put $d = r+s-1 (=\dim A)$ and set  
$$
\frm =(x_1,\,\ldots,x_r)R,\; \frn = (y_1,\ldots,y_s)S,\quad \text{and} \;\;
\fraM = \frm \# \frn = \bigoplus_{n=1}^{\infty} R_n \otimes S_n. 
$$
Then the (graded) canonical module $K_A$ of $A$ is isomorphic to $K_R \#K_S$ 
by \cite{GW, Theorem 4.3.1}
(In particular, $A$ is Gorenstein if and only if $r=s$). 
Thus by virtue of \cite{GW,Theorem 4.1.5}, we get  
$$
E_A = H_{\fraM}^d(K_A) = H_{\fraM}^d(K_R \#K_S) 
= H_{\frm}^r(K_R) \# H_{\frn}^s(K_S) 
= E_R \# E_S.  
$$
\par 
Further, since $E_R$ (resp. $E_S$) can be represented as a graded module 
$k[x_1^{-1},\,\ldots,x_r^{-1}]$ (resp. $k[y_1^{-1},\,\ldots,y_s^{-1}]$)
which is called {\it the inverse system of Macaulay}, we get 
$$
E_A \cong k[x_1^{-1},\,\ldots,x_r^{-1}] \# k[y_1^{-1},\,\ldots,y_s^{-1}]. 
$$
Then $z= 1\# 1 \in E_A$ generates the socle of $E_A$. 
\medskip 
Using this, we have 
\proclaim {Proposition 5.5}
Let $A=R\#S$ and $z = 1 \#1$ be as above. 
Then
$$
\multlinegap{1.3cm}
\multline 
  l_A(A/\ann_A(F^e_A(z)) \\
= \#\left\{(a_1,\dots,a_r,b_1,\dots,b_s) 
\in \bbZ^{r+s} \;\Biggm|\;
\aligned 
&0 \le a_1,\,\dots,a_r \le q-1 \\
&0 \le b_1,\,\dots,b_s  \le q-1 \\
&a_1 + \dots + a_r  = b_1 + \dots + b_s  
\endaligned 
\right\}. 
\endmultline
$$
\endproclaim 
\demo{\quad \bf Proof}
We use the same notation as in the above argument. 
Now we shall investigate the Frobenius action of $z$ in $E_A$. 
First note that $\FA^e(E_A) \cong \FR^e(E_R) \# \FS^e(E_S)$. 
Thus it is enough to investigate the Frobenius action of $z_1=1$ in $E_R$. 
Since $E_R = H_{\frm}^r(R)(-r)$, that is, 
$H_{\frm}^r(R) \cong (x_1 \cdots x_r)^{-1}E_R$, 
the generator $z_1$ of $E_R$ corresponds to the element 
$w_1= (x_1 \cdots x_r)^{-1}$ via this isomorphism. 
Then we have $F_R^e(w_1) = (x_1  \cdots x_r)^{-q}$ 
since there exsits an isomorphism 
$$
\align
&(x_1 \cdots x_r)^{-1}  k[x_1,\,\ldots,x_r] \to 
H_{\frm}^{r}(R) = \varinjlim_{n} R/(x_1^n,\,\ldots,x_d^n). \\
& \left(x_1^{-a_1}\cdots x_d^{-a_d} \mapsto 
[x_1^{a-a_1}\cdots x_d^{a-a_d} + (\underline{x}^a)], 
\text{where}\;\; a := \max\{a_1,\,\ldots,a_d\}\right)
\endalign 
$$
If we identify $\FR^e(E_R)$ with $E_R$,  
then $$
F_R^e(z_1) = (x_1\cdots x_r) \cdot F^e(w_1) 
= (x_1 \cdots x_r)^{-(q-1)}.
$$ 
Therefore 
$$
F_A^e(z) = F_R^e(z_1) \#  F_S^e(z_2) 
= (x_1 \cdots x_r)^{-(q-1)} \# (y_1 \cdots y_s)^{-(q-1)} 
\quad \text{in $E_A$}. 
\tag"\bf (5.5.1)"
$$
\par
For any element $c = x_1^{a_1}\cdots x_r^{a_r} \# y_1^{b_1}\cdots 
y_s^{b_s}$ in $R$, we have 
$$  
cF^e(z) \ne 0 \quad \text{in}\; E_A \qquad 
\Longleftrightarrow \qquad 
\left\{
\aligned 
&0 \le a_1,\,\dots,a_r \le q-1, \\
&0 \le b_1,\,\dots,b_s  \le q-1, \\
&a_1 + \dots + a_r  = b_1 + \dots + b_s.  
\endaligned 
\right.
$$
Thus we get the required assertion. 
\qed 
\enddemo
\medskip 
We are now ready to state our main theorem in this section. 
\proclaim {Theorem 5.6}
Let $A=k[x_1,\,\ldots,x_r]\#k[y_1,\,\ldots,y_s]$, where 
$2 \le r \le s$, and put $d = r+s-1$. 
Then 
$$ 
\mhk(A)  =\frac{r!}{d!}\,S(d,r) 
+\frac{1}{d!}\sum_{k=1}^{r-1} \sum_{j=1}^{r-k}
\binom{r}{k+j} \binom{s}{j}(-1)^{r+k}k^d, 
$$
where $S(n,k)$  denotes the Stirling number of the second kind$;$ see below. 
\par 
In particular, 
$$
\ehk(A) + \mhk(A) = \dfrac{r! \cdot S(d,r)+s! \cdot S(d,s)}{d!}. 
$$
\endproclaim 
\par 
The following two corollaries easily follow from Theorem 5.1 and Theorem 5.6. 

\proclaim {Corollary 5.7}
Let $A = R \# S = k[x_1,\,x_2]\# k[y_1,\ldots,y_s]$, which is isomorphic 
to the Rees algebra $S[\frn t]$ over $S$. 
Then 
$$
\mhk(A) = \frac{2^{s+1}-s-2}{(s+1)!}. 
$$  
\endproclaim 

\proclaim {Corollary 5.8}
Under the same notation as in Thereom 5.6, 
further, assume that $A$ is Gorenstein, that is, $r=s$. 
Then 
$$
\ehk(A) + \mhk(A) = \dfrac{2 \cdot r!}{(2r-1)!}S(2r-1,r). 
$$  
\endproclaim 
\demo{\quad \bf Proof of Theorem 5.6}
If we put 
$\alpha_{r,n} := l_R(\frm^n/\frm^{n+1}) = \binom{n+r-1}{r-1}$ 
and 
$\alpha_{r,n,q} := l_R(\frm^n/\frm^{n-q}\frm^{[q]}+\frm^{n+1})$, 
then 
$$
\align
\ehk(A) & = 
\lim_{q \to \infty} 
\dfrac{1}{q^d}\sum_{n=0}^{r(q-1)}\alpha_{r,n}\alpha_{s,n,q} \\
& + \lim_{q \to \infty} 
\dfrac{1}{q^d}\sum_{n=0}^{s(q-1)}\alpha_{r,n,q}\alpha_{s,n}
- \lim_{q \to \infty} 
\dfrac{1}{q^d}\sum_{n=0}^{r(q-1)}\alpha_{r,n,q}\alpha_{s,n,q}.
\endalign 
$$
Also, by virtue of Proposition 5.5, we get 
$$
\mhk(A) = \lim_{q \to \infty} 
\dfrac{1}{q^d}\sum_{n=0}^{r(q-1)}\alpha_{r,n,q}\alpha_{s,n,q}. 
$$
Hence the required assertion follows from the following lemma. \qed 
\enddemo 

\proclaim {Lemma 5.9} {\rm (cf. \cite{EY, Lemma 3.8, Lemma 3.9})}
Under the same notation as above, we have 
$$
\align 
(1) \quad \lim_{q \to\infty}\frac{1}{q^d}\sum_{n=0}^{r(q-1)}
\alpha_{r,n,q}\alpha_{s,n} 
& =\frac{r!}{d!}S(d,r).\\ 
& \\ 
(2) \quad \lim_{q \to\infty}\frac{1}{q^d}
\sum_{n=0}^{r(q-1)}
\alpha_{r, n, q}\alpha_{s, n, q} 
& =\frac{r!}{d!} S(d,r)  +\frac{1}{d!}\sum_{0<j<i\leq r}
\binom{r}{i} \binom{s}{j}(-1)^{r-i+j}(i-j)^{d}.
\endalign 
$$
\endproclaim 

%

\enddocument